\documentclass[11pt,reqno]{amsart}
\usepackage{amsthm}
\usepackage{amssymb}
\usepackage{latexsym}
\usepackage{multicol}
\usepackage{verbatim,enumerate}
\usepackage{amscd}
\usepackage{hyperref}
\usepackage{amsmath, amscd}

\advance\textwidth by 1.2in \advance\oddsidemargin by -.6in \advance\evensidemargin by -.6in
\newtheorem*{lem}{Lemma}
\newtheorem*{prop}{Proposition}

\theoremstyle{definition}

\theoremstyle{definition}
\newtheorem{thm}{Theorem}

\newenvironment{pf}{\proof}{\endproof}
\newcounter{cnt}
\newenvironment{enumerit}{\begin{list}{{\hfill\rm(\roman{cnt})\hfill}}{%
\settowidth{\labelwidth}{{\rm(iv)}}\leftmargin=\labelwidth%
\advance\leftmargin by \labelsep\rightmargin=0pt\usecounter{cnt}}}{\end{list}} \makeatletter
\def\mydggeometry{\makeatletter\dg@YGRID=1\dg@XGRID=20\unitlength=0.003pt\makeatother}
\makeatother \theoremstyle{remark}

\numberwithin{equation}{section}

\let\bwdg\bigwedge
\def\bigwedge{{\textstyle\bwdg}}

\begin{document}

\newcommand{\thmref}[1]{Theorem~\ref{#1}}
\newcommand{\secref}[1]{Section~\ref{#1}}
\newcommand{\lemref}[1]{Lemma~\ref{#1}}
\newcommand{\propref}[1]{Proposition~\ref{#1}}
\newcommand{\corref}[1]{Corollary~\ref{#1}}
\newcommand{\remref}[1]{Remark~\ref{#1}}
\newcommand{\defref}[1]{Definition~\ref{#1}}
\newcommand{\er}[1]{(\ref{#1})}
\newcommand{\id}{\operatorname{id}}
\newcommand{\ord}{\operatorname{\emph{ord}}}
\newcommand{\sgn}{\operatorname{sgn}}
\newcommand{\wt}{\operatorname{wt}}
\newcommand{\tensor}{\otimes}
\newcommand{\from}{\leftarrow}
\newcommand{\nc}{\newcommand}
\newcommand{\rnc}{\renewcommand}
\newcommand{\dist}{\operatorname{dist}}
\newcommand{\qbinom}[2]{\genfrac[]{0pt}0{#1}{#2}}
\nc{\cal}{\mathcal} \nc{\goth}{\mathfrak} \rnc{\bold}{\mathbf}
\renewcommand{\frak}{\mathfrak}
\newcommand{\supp}{\operatorname{supp}}
\newcommand{\Irr}{\operatorname{Irr}}
\renewcommand{\Bbb}{\mathbb}
\nc\bomega{{\mbox{\boldmath $\omega$}}} \nc\bpsi{{\mbox{\boldmath $\Psi$}}}
 \nc\balpha{{\mbox{\boldmath $\alpha$}}}
 \nc\bpi{{\mbox{\boldmath $\pi$}}}
 \nc\bxi{{\mbox{\boldmath $\xi$}}}
  \nc\bmu{{\mbox{\boldmath $\mu$}}}
\nc\bsigma{{\mbox{\boldmath $\sigma$}}} \nc\bcN{{\mbox{\boldmath $\cal{N}$}}} \nc\bcm{{\mbox{\boldmath $\cal{M}$}}} \nc\bLambda{{\mbox{\boldmath
$\Lambda$}}}

\newcommand{\lie}[1]{\mathfrak{#1}}
\makeatletter
\def\section{\def\@secnumfont{\mdseries}\@startsection{section}{1}%
  \z@{.7\linespacing\@plus\linespacing}{.5\linespacing}%
  {\normalfont\scshape\centering}}
\def\subsection{\def\@secnumfont{\bfseries}\@startsection{subsection}{2}%
  {\parindent}{.5\linespacing\@plus.7\linespacing}{-.5em}%
  {\normalfont\bfseries}}
\makeatother
\def\subl#1{\subsection{}\label{#1}}
 \nc{\Hom}{\operatorname{Hom}}
  \nc{\mode}{\operatorname{mod}}
\nc{\End}{\operatorname{End}} \nc{\wh}[1]{\widehat{#1}} \nc{\Ext}{\operatorname{Ext}} \nc{\ch}{\text{ch}} \nc{\ev}{\operatorname{ev}}
\nc{\Ob}{\operatorname{Ob}} \nc{\soc}{\operatorname{soc}} \nc{\rad}{\operatorname{rad}} \nc{\head}{\operatorname{head}}
\def\Im{\operatorname{Im}}
\def\gr{\operatorname{gr}}
\def\mult{\operatorname{mult}}
\def\Max{\operatorname{Max}}
\def\ann{\operatorname{Ann}}
\def\sym{\operatorname{sym}}
\def\Res{\operatorname{\br^\lambda_A}}
\def\und{\underline}
\def\Lietg{$A_k(\lie{g})(\bsigma,r)$}

 \nc{\Cal}{\cal} \nc{\Xp}[1]{X^+(#1)} \nc{\Xm}[1]{X^-(#1)}
\nc{\on}{\operatorname} \nc{\Z}{{\bold Z}} \nc{\J}{{\cal J}} \nc{\C}{{\bold C}} \nc{\Q}{{\bold Q}}
\renewcommand{\P}{{\cal P}}
\nc{\N}{{\Bbb N}} \nc\boa{\bold a} \nc\bob{\bold b} \nc\boc{\bold c} \nc\bod{\bold d} \nc\boe{\bold e} \nc\bof{\bold f} \nc\bog{\bold g}
\nc\boh{\bold h} \nc\boi{\bold i} \nc\boj{\bold j} \nc\bok{\bold k} \nc\bol{\bold l} \nc\bom{\bold m} \nc\bon{\bold n} \nc\boo{\bold o}
\nc\bop{\bold p} \nc\boq{\bold q} \nc\bor{\bold r} \nc\bos{\bold s} \nc\boT{\bold t} \nc\boF{\bold F} \nc\bou{\bold u} \nc\bov{\bold v}
\nc\bow{\bold w} \nc\boz{\bold z} \nc\boy{\bold y} \nc\ba{\bold A} \nc\bb{\bold B} \nc\bc{\bold C} \nc\bd{\bold D} \nc\be{\bold E} \nc\bg{\bold
G} \nc\bh{\bold H} \nc\bi{\bold I} \nc\bj{\bold J} \nc\bk{\bold K} \nc\bl{\bold L} \nc\bm{\bold M} \nc\bn{\bold N} \nc\bo{\bold O} \nc\bp{\bold
P} \nc\bq{\bold Q} \nc\br{\bold R} \nc\bs{\bold S} \nc\bt{\bold T} \nc\bu{\bold U} \nc\bv{\bold V} \nc\bw{\bold W} \nc\bz{\bold Z} \nc\bx{\bold
x} \nc\KR{\bold{KR}} \nc\rk{\bold{rk}} \nc\het{\text{ht }}

\nc\toa{\tilde a} \nc\tob{\tilde b} \nc\toc{\tilde c} \nc\tod{\tilde d} \nc\toe{\tilde e} \nc\tof{\tilde f} \nc\tog{\tilde g} \nc\toh{\tilde h}
\nc\toi{\tilde i} \nc\toj{\tilde j} \nc\tok{\tilde k} \nc\tol{\tilde l} \nc\tom{\tilde m} \nc\ton{\tilde n} \nc\too{\tilde o} \nc\toq{\tilde q}
\nc\tor{\tilde r} \nc\tos{\tilde s} \nc\toT{\tilde t} \nc\tou{\tilde u} \nc\tov{\tilde v} \nc\tow{\tilde w} \nc\toz{\tilde z} \nc\woi{w_{\omega_i}}

\title{An application of global Weyl modules of $\lie{sl}_{n+1}[t]$ to  invariant theory.}
\author{Vyjayanthi Chari and Sergei Loktev}
\thanks{\scriptsize{VC acknowledges support from the NSF, DMS-0901253}.\\ {\scriptsize SL is supported by grants RFBR-09-01-00239, RFBR-CNRS-09-02-93106, N.Sh.8462.2010.1
} }
\email{chari@math.ucr.edu}
\email{s.loktev@gmail.com}

\maketitle
\begin{abstract} We identify $\lie{sl}_{n+1}$ isotypical components of global Weyl
modules with natural subspaces in a polynomial ring, and then apply the 
representation theory of current algebras to classical problems in invariant theory.

\end{abstract}
\section*{Introduction}
Given  a positive integer  $k$, let $A_{k}$ be
the algebra of polynomials with complex coefficients in $k$
variables $t_1,\cdots, t_{k}$ and let $S_{k}$ be the symmetric group
on $k$ letters. The group $S_{k}$ acts naturally on $ A_{k}$ by
permuting the variables. It is a classical result that  $A_k$ when regarded as a module for (the algebra of invariants)
$A_k^{S_k}$ is free of rank $k!$.

In this paper, we are interested in certain variations of this
problem. We regard $A_k$ as being graded  by the non--negative integers in the standard way with each generator  having grade one and note that $A_k^{S_k}$ is a graded subalgebra. Given a partition  $\xi=\xi_1\ge\xi_2\ge\cdots\ge \xi_{n+1}>0$ of $k$  we define polynomials $\bop(\bor)\in A_k$ indexed by elements in $\bz_+^{\xi_2}\times\cdots\times\bz_+^{\xi_{n+1}}$, where $\bz_+$ is the set of non--negative integers. (These polynomials can be thought of as being  associated to Young symmetrizers of tableaux of shape  $\xi$ where the filling in the first row  is $1$).
For instance, when  $\xi=k-\ell\ge\ell >0$, the polynomials are indexed by elements of $\bz_+^\ell$ and are defined by,
\begin{equation}\label{defnp}
\bop(\bos)=\sum_{\sigma\in
S_\ell}\left(\prod_{i=1}^\ell(t_{2i-1}^{s_{\sigma(i)}}-t_{2i}^{s_{\sigma(i)}})\right),\ \ \bos\in\bz_+^\ell.
\end{equation}
Let $M_{k,\xi}$ be the $A_k^{S_k}$--submodule of $A_k$ generated by the elements $\bop(\bor)$ $\bz_+^{\xi_2}\times\cdots\times\bz_+^{\xi_{n+1}}$. 
We use  the theory  global Weyl modules of the current algebra of $\lie{sl}_{n+1}$ to prove that $M_{k,\xi}$ is free and has a graded basis
where the number of elements of a given grade is the coefficient of the corresponding power of $q$ in the Kostka polynomial $\kappa_{1^k, \xi^{tr}}$ and show that it can be identified with the 
 multiplicity space of $\xi$ in $A_k$.
  
 The global Weyl modules we use are indexed by  multiples of the first fundamental weight of $\lie{sl}_{n+1}$. These modules are well--understood in the literature and many spanning sets and bases are known for them: a Poincare--Birkhoff-Witt type spanning set, a Gelfand--Tsetlin basis (\cite{CL}) a global canonical bases (\cite{K}. In fact our main result  gives an explicit way to construct a $A_k^{S_k}$--generating set or basis of $M_{k,\xi}$ from a spanning set or basis of the appropriate global Weyl module.
The special case described above corresponds to taking a Poincare--Birkhoff--Witt spanning set for the global Weyl modules.

For $\xi=\ell\ge\ell$, the result was
first conjectured by Chari and Greenstein. It was motivated by their
study  \cite{CG1}, \cite{CG2} of the homological properties of the
category of finite--dimensional representations of  the current algebra of $\lie{sl}_2$. In \cite{BCDM} the module $M_{k,\ k-\ell\ge\ell}$ is denoted as
$\bm_{k,\ell}$ and a more specific conjecture is given which
identifies a natural  basis for it and this conjecture remains open.

Our interest in this result stems from the fact that it plays a crucial role in \cite{BCM}
which establishes a BGG--type duality for locally
finite--dimensional representations of the current algebra associated to $\lie{sl}_2$. In Section 1, we give an elementary formulation of the results of this paper  without using the notion of global Weyl modules.  The definitions and necessary  results on global Weyl modules are  recalled in  Section 2.  In  Section 3 we prove our main  result which makes the identification between the module $M_{k,\xi}$ and the appropriate space of invariants in the global Weyl modules.

{\em Acknowledgements.} We thank Shrawan Kumar for several discussions.  This work was completed when both authors were visiting the  the Hausdorff Institute, Bonn, in connection with the trimester \lq\lq On the interactions of Representation theory with Geometry and Combinatorics,\rq\rq
and the authors   thank the organizers of the trimester for the invitation to participate in the program.

\section{The Modules $M_{k,\xi}$}

\subsection{} Throughout this note we let $\bz_+$ (resp. $\bn$)  denote the set of non--negative (resp. positive)  integers and $\bc$ the set of complex numbers. For $k\in\bn$, we let $A_k$ be  the polynomial ring in $k$--variables equipped with the standard grading. For $r\in\bz_+$ the $r$--th graded piece is denoted as  $A_k[r]$ and  is the space of homogenous polynomials of degree $r$. We shall use without further mention, the celebrated result of Quillen (see \cite{L} for an exposition) that any projective module for a polynomial ring is free.

\subsection{}\label{gp} We begin with some standard  results \cite{GW} on the  representation theory of the symmetric group $S_k$ on $k$ letters. Let $\cal P[k]$ be the set of all partitions $\xi=\xi_1\ge\cdots\ge\xi_{n+1}\ge 0$  of $k$.
The irreducible representations  of $S_k$ are parametrized by elements of $\cal P[k]$.
 By abuse of notation, we shall denote by  $\xi$ both a partition of $k$ and an irreducible representation of $S_k$ in the corresponding isomorphism class and we set $d_\xi=\dim\xi$.
 Any representation $V$ of $S_k$  can be written as a sum of finite dimensional irreducible modules and we  denote by $V^\xi$  the isotypical component corresponding to $\xi\in\cal P[k]$, and we have    $V=\oplus_{\xi\in\cal P[k]} V^\xi$ as $S_k$--modules. We shall also use then notation $V^{S_r}$ to denote the isotypical component corresponding to the trivial representation.
 Let $\bc[S_k]$ be  the group algebra of $S_k$ over $\bc$ and  regard it as module for itself via right multiplication. Then $$\bc[S_k]=\bigoplus_{\xi\vdash k}\bc[S_k]^\xi, \ \ \  \ \dim\Hom_{S_k}(\xi,\bc[S_k]^\xi)=d_\xi.$$
  The following is well--known but we isolate is as a Lemma since it is used frequently. Recall that if $V$ and $V'$ are two representations of $S_k$ then one has the diagonal action of $S_k$ on $V\otimes V'$.
\begin{lem} We have, $\Hom_{S_k}(\xi,\xi)\ne 0$ or equivalently that $(\xi\otimes\xi)^{S_k}\ne 0$. Suppose that   $V$ is any representation of $S_k$  such that $V^\xi\ne 0$. Then any non-zero element $v\in (V^\xi\otimes \xi)^{S_k}$ can be written as a sum, $$v=\sum_{s=1}^{d_\xi}w_s\otimes v_s, w_s\in V^\xi,\ \ v_s\in \xi,$$ where the elements $\{v_s:1\le s\le d_\xi\}$ are a basis of $\xi$ and  the elements $\{w_s$: $1\le s\le d_\xi\}$ span an irreducible $S_k$--submodule of $V^\xi$ (in particular the elements $\{w_s$, $1\le s\le \xi\}$ are linearly independent). \hfill\qedsymbol
\end{lem} 

\subsection{} Consider the natural right action of $S_k$ on the ring $A_k$ given by permuting the variables. The action preserves the graded pieces of $A_k$ and the subalgebra $A_k^{S_k}$ of invariants is a graded subalgebra of $A_k$. Regard $A_k$ as a module over the ring $A_k^{S_k}$ via right multiplication and note that it commutes with the action of $S_k$.  Hence for any partition $\xi$  of $k$,  the $S_k$--isotypical component $A_k^\xi$ is also a module for $A_k^{S_k}$ and we have an isomorphism of $A_k^{S_k}$ and $S_k$--modules,
$$A_k\cong \bigoplus_{\xi\vdash k} A_k^\xi.$$ By a classical result of Chevalley one knows that $A_k^{S_k}$ is also a polynomial ring in $k$ variables and if we denote  by $I_k$ the augmentation ideal in $A_k^{S_k}$, then one has an isomorphism of  $S_k$--modules, $$A_k/I_kA_k^{S_k}\cong \bc[S_k].$$
 It follows that $A_k$ is a free module for $A_k^{S_k}$ of rank $k!$ and in fact  that each isotypical component $A_k^\xi$ is also a free $A_k^{S_k}$--module of rank $d_\xi^2$. The augmentation ideal $I_k$ is graded and hence $\bc[S_k]$ acquires a $\bz_+$--grading. It is known that the graded multiplicity of the irreducible representation $\xi$ is given by the Kostka polynomial $\kappa_{1^{tr},\xi}$. 

\subsection{}\label{gen} To state our first result fix a partition $\xi=\xi_1\ge\cdots\ge\xi_{n=1}>0$ of $k$ with  $(n+1)$ parts and set $k_i=\xi_i-\xi_{i+1}$, $1\le i\le n+1$ where we understand that $k_{n+1}=\xi_{n=1}$. We shall define a family of elements $\bop(r)$ in $A_k$ which are indexed by elements of the set $\bz_+^{\xi_2}\times\cdots\times\bz_+^{\xi_{n+1}}$.
It  will be convenient to think of  $A_k$ as a tensor product of algebras  $$A_k=A_1^{\otimes k_1}\otimes A_2^{\otimes k_2}\otimes \cdots\otimes A_{n+1}^{\otimes k_{n+1}},$$ and we shall also identify
$$\Phi_\xi: \bz_+^{\xi_2}\times\cdots\times\bz_+^{\xi_{n+1}}\cong \bz_+^{k_2}\times(\bz^2_+)^{k_3}\cdots\times (\bz_+^{n+1})^{k_{n+1}},$$ as follows. (We use the convention that $\bz_+^0$ is the emptyset.) An element on the left hand side can naturally be thought of as an upper triangular matrix of size $k\times\xi_2$ (the entries in the $(r,s)$--the place are defined to be zero if  $1\le s\le k_r$) and this maps to the transpose matrix and then we drop all the zero entries. For example, if $\xi= 3\ge 3\ge 1$, then the map is  $$((r_1,r_2,r_3), s_1) \longrightarrow (r_1,r_2, (r_3,s_1)),$$  while if $\xi=2\ge 2\ge 2$, then the map is $$((r_1,r_2),(s_1,s_2))\to ((r_1,s_1), (r_2,s_2)).$$
For $2\le \ell \le k$ and $1\le s\le k_\ell$, let $\Phi_\xi^{\ell,s}: \bz_+^{\xi_2}\times\cdots\times\bz_+^{\xi_{n+1}}\to(\bz_+^\ell)^{k_\ell}\to \bz_+^\ell$ be the projection onto the $\ell$--th factor and then onto the $s$--th factor. Given $(r_1,\cdots, r_\ell)\in\bz_+^\ell$ set $$\boa(r_1,\cdots, r_\ell)=\sum_{\tau\in S_{\ell+1}}{\rm{sgn}}(\tau) t^{r_1}_{\tau(1)}\cdots t_{\tau(\ell)}^{r_\ell}\in A_{\ell+1}.$$
Next, given  $\bor\in\bz_+^{\xi_2}\times\cdots\times\bz_+^{\xi_{n+1}}$, set
  $$\bob(\bor)=1^{\otimes k_1}\otimes \boa(\Phi_\xi^{2,1}(\bor))\otimes\cdots\otimes\boa(\Phi_\xi^{2,k_2}(\bor))\otimes\cdots\otimes \boa(\Phi_\xi^{n+1,1}(\bor))\otimes
\cdots\otimes\boa(\Phi^{n+1,k_{n+1}}_\xi(\bor)).$$ Finally, we set $$\bop(\bor)=\sum_{\sigma\in\bs_\xi}\bob(\sigma(\bor)),$$where $\bs_\xi$ is the group $S_{\xi_2}\times\cdots\times S_{\xi_{n+1}}$ acting naturally on $\bz_+^{\xi_2}\times\cdots\times\bz_+^{\xi_{n+1}}$.

 \subsection{} We compute  an example  to illustrate the preceding definitions and to relate it to the modules in the introduction.   Suppose that $n=1$ and  $\xi=k-m\ge m>0$ so that  $k_1=k-2m$, $k_2=m$. The polynomials we are interested in are indexed by elements of $\bor=(r_1,\cdots ,r_m)\in\bz_+^m$ and the polynomials live in $A_1^{\otimes k-m}\otimes A_2^{\otimes m}$. Since $\Phi_\xi^{2,\ell}(\bor)=r_\ell$,we have  $$\bob(\bor)=1^{k-m}\otimes \boa(r_1)\otimes\cdots\otimes\boa(r_m)\ \qquad \boa(r)=(t_1^r-t_2^r)\in A_2.$$ Hence we get $$\bop(\bor)=1^{\otimes (k-m)}\otimes \sum_{\sigma\in S_m}(t_1^{r_{\sigma(1)}}-t_2^{\sigma){r_\sigma(1)}})\cdots (t_1^{r_{\sigma(m)}}-t_2^{\sigma){r_\sigma(m)}}).$$

\subsection{} We shall prove,
\begin{thm} Let $k\in\bn$ and let $\xi=\xi_1\ge\xi_2\ge\cdots\ge\xi_{n+1}>0$ be a partition of $k$ into  $(n+1)$ parts and let $k_i=\xi_i-\xi_{i+1}$. The  $A_k^{S_k}$--submodule of $A_k$ (denoted $M_{k,\xi})$ generated by the elements $\bop(\bor)$, $\bor\in \bz_+^{\xi_2}\times\cdots\times\bz_+^{\xi_{n+1}}$,is  graded free   and the number of elements of grade $s$ in the basis is given the coefficient of $q^s$ in the Kostka polynomial $\kappa_{1^k,\xi^{tr}}$.
 Moreover, we have an isomorphism of graded $S_k$ and $A_k^{S_k}$--modules $$A_k^\xi\cong \cong M_{k,\xi}\otimes \xi,$$ where we regard $A_k^{S_k}$ on $M_{k,\xi}$ and $S_k$ as acting on $\xi$.  In particular, if $I_k$ is the augmentation ideal in $A_k$, we have a $\bz_+$--graded map $$M_{k,\xi}/I_kM_{k,\xi}\cong \Hom_{S_k}(\xi,\bc[S_k]).$$

\end{thm}
\medskip
We deduce this theorem from the results of the next two sections .
\section{Representations of $\lie{sl}_{n+1}[t]$}

In this section we first recall several standard results on the representation theory of $\lie{sl}_{n+1}$ which can be found in any standard book \cite{GW}, \cite{H}  for instance. We then  define a family of representations (the global Weyl modules) of the infinite--dimensional Lie algebra $\lie{sl}_{n+1}\otimes\bc[t]$ and recall some of their properties which were  first established in \cite{CPweyl} and  \cite{FL}.

 \subl{}\label{slrep} Let  $\lie {sl}_{n+1}$  be the Lie algebra of $(n+1)\times (n+1)$--matrices matrices of trace zero and  $\lie h$ (resp. $\lie n^+$, $\lie n^-$) be the subalgebra of diagonal (resp.  strictly upper, lower triangular) matrices. Let $e_{i,j}$ be the $(n+1)\times (n+1)$--matrix with one in the $(i,j)$--th position and zero elsewhere. The elements $$x_i^+= e_{i,i+1},\ \ h_i=e_{i,i}-e_{i+1,i+1},\ \ x_i^-= e_{i+1,i},\ \  1\le i\le n,$$ generate $\lie{sl}_{n+1}$ as a Lie algebra. Let $\Omega:\lie{sl}_{n+1}\to\lie{sl}_{n+1}$ be the anti--automorphism satisfying \begin{equation}\label{chevant} \Omega(x_i^\pm)=x^\mp_i,\ \ \ \ \ \ \Omega(h_i)=h_i.\end{equation}

\subsection{} Let  $\cal P^n$ be  the set of all partitions with at most $(n+1)$--parts and for $1\le i\le n$, let $\omega_i\in\cal P^n$ be the partition $1^i$. An element  $\lambda=\lambda_1\ge\lambda_2\ge\cdots\ge \lambda_{n+1}\ge 0$ of $\cal P^n$ defines a finite--dimensional irreducible representation $V(\lambda)$ of $\lie{sl}_{n+1}$ and  two partitions $\lambda,\lambda'$ define isomorphic representations iff  $\lambda_i-\lambda_{i+1}=\lambda_i'-\lambda_{i+1}'$ for all $ 1\le
 i\le n$. The trivial representation of $\lie{sl}_{n+1}$ will be denoted as $V(\omega_0)$ and corresponds to taking the empty partition.

  The module $V(\lambda)$ is generated by an element $v_\lambda$ with defining relations:
$$\lie n^+ v_\lambda=0,\ \ \ h_i v_\lambda=(\lambda_i-\lambda_{i+1})v_\lambda,\ \ \ (x_i^-)^{\lambda_i-\lambda_{i+1}+1}v_\lambda=0,\ \ \ \ 1\le i\le n.$$
 Any irreducible finite--dimensional representation of $\lie{sl}_{n+1}$ is isomorphic to $V(\lambda)$ for some $\lambda\in\cal P^n$.

  \subsection{} We say that a representation $V$ of $\lie{sl}_{n+1}$ is locally finite--dimensional if it is a direct sum of  finite--dimensional representations. By Weyl's theorem, this is equivalent to requiring that $V$ is isomorphic to a direct sum of $V(\lambda)$, $\lambda\in\cal P^n$.
 The Lie algebra $\lie h$ acts diagonalizably on  a locally finite--dimensional module $V$ with integer eigenvalues,
 and we have $$V=\bigoplus_{\bor\in\bz^n} V_\bor, \quad \ V_\bor=\{v\in V: h_i v= r_i v\},\quad  \ \bor\in\bz^n.$$ In the case when $\xi$ is a partition, we also use the notation,
 \begin{gather*}V_\xi=\{v\in V: \ \ h_i v=(\xi_i-\xi_{i+1})v,\ \ 1\le i\le n\},\qquad \ V_\xi^{\lie n^+}=\{v\in V_\xi:\lie n^+ v=0\}.\end{gather*} Note that $V(\lambda)^{\lie n^+}=\bc v_\lambda$ for all $\lambda\in\cal P^n$ and that in general, a locally finite--dimensional representation is generated as a $\lie{sl}_{n+1}$--module by the spaces $V_\xi^{\lie n^+}$, $\xi\in\cal P^n$. The following is elementary.
 \begin{lem} \label{ninv}Let $V$ be a locally finite--dimensional $\lie{sl}_{n+1}$--module and let $\xi\in\cal P^n$. Then $$V_\xi= V_\xi^{\lie n^+}\oplus (\lie n^-V \cap V_\xi).$$\hfill\qedsymbol
 \end{lem}

\subsection{}  Recall that if $U$ is any vector space there is a natural action of $S_k$ on the $k$--fold tensor product $U^{\otimes k}$ which just permutes the factors. We shall need the following well--known  result, the first two parts are elementary and the third part is the famous Schur--Weyl duality between representations of $\lie{sl}_{n+1}$ and those of $S_k$

\begin{prop}\label{schweyl} Let $k,n\in\bn$.
\begin{enumerit}
\item[(i)] The module $V(\omega_1)$ is isomorphic to the  natural  representation $\bc^{n+1}$ of $\lie{sl}_{n+1}$ with standard basis $e_1,\cdots, e_{n+1}$. Moreover the assignment $(e_i,e_j)=\delta_{i,j}$ defines a symmetric bilinear from $(\ ,\ ): V(\omega_1)\otimes V(\omega_1)\to\bc$ satisfying,
    $$(xv, v')=(v, \Omega(x)v'),\ \ v,v'\in V(\omega_1),\ \ \ x\in\lie{sl}_{n+1}.$$
 \item[(ii)] For $1\le i\le n+1$, we have $$(V(\omega_1)^{\otimes i}))^{\lie n^+}_{\omega_i}=\bc (e_1\wedge e_2\wedge\cdots\wedge e_i),$$ and an isomorphism of $\lie{sl}_{n+1}$--modules $\bigwedge^i V(\omega_1)\cong V(\omega_i)$.
     \item[(iii)] Let $k\in\bn$. The natural right action of $S_k$ on $V(\omega_1)^{\otimes k}$ commutes with the left action of $\lie{sl}_{n+1}$ and
 as a $(\lie{sl}_{n+1}, S_k)$--bimodule, we have,
     $$V(\omega_1)^{\otimes k}\cong\bigoplus_{\xi\in\cal P^{n}\cap\cal P[k]}V(\xi)\otimes \xi,$$ where the sum is over all partitions of $k$ with at most $(n+1)$ parts. Equivalently, the space $(V(\omega_1)^{\otimes k})^{\lie n^+}_\xi$ is an irreducible $S_k$--submodule of $V(\omega_1)^{\otimes k}$ and  we have an isomorphism of $S_k$--modules $$(V(\omega_1)^{\otimes k})^{\lie n^+}_\xi\cong \xi.$$
\end{enumerit} \hfill\qedsymbol\end{prop}

\subsection{}\label{form1} Denote by $(\ ,\ )_{\scriptsize k}:V(\omega_1)^{\otimes k}\times V(\omega_1)^{\otimes k}\to\bc$  the symmetric non--degenerate bilinear form given by extending linearly the assignment, $$(e_{i_1}\otimes \cdots\otimes e_{i_k},\   e_{j_1}\otimes\cdots\otimes e_{j_k})_{\scriptsize k}= (e_{i_1}, e_{j_1})\cdots (e_{i_k}, e_{j_k}),$$ where $ i_s, j_s\in\{1,\cdots, n+1\}$ for $1\le s\le k$.  The form is clearly $S_k$--invariant and moreover, $$(x\bov,\  \bov')_{\scriptsize k} = (\bov,\  \Omega(x) \bov')_{\scriptsize k},\quad \ \bov,\bov'\in V(\omega_1)^{\otimes k}, \quad  \ x\in\lie g.$$ In particular since $\Omega(h_i)=h_i$ for all $1\le i\le n$, we have $$(V(\omega_1)^{\otimes k}_\bor, \ V(\omega_1)^{\otimes k}_{\bor'})_{\scriptsize k}=0,\ \ \bor,\bor'\in\bz^n, \ \ \bor\ne\bor',$$ and so the restriction of the form to $V(\omega_1)^{\otimes k}_\bor$ is non--degenerate.
\begin{lem}\label{nondeg} For all  $\xi\in\cal P^n\cap\cal P[k]$ the restriction of the form $(\ ,\ )_{\scriptsize k}$ to $(V(\omega_1)^{\otimes k})^{\lie n^+}_\xi\times(V(\omega_1)^{\otimes k})^{\lie n^+}_\xi$ is non--degenerate.
\end{lem}
\begin{pf} Let  $\bov,\bov'\in V(\omega_1)^{\otimes k}_\xi$ be such that $(\bov,\ \bov')_{\scriptsize k}\ne 0$ and assume that $\lie n^+\bov=0$. By Lemma \ref{ninv}, we may write $\bov'= \bov_1 +\bov_2,$ where $\lie n^+\bov_1=0$ and $\bov_2 =\sum_{i=1}^ne_{i+1,i}\bou_i$ for some $\bou_i\in  V(\omega_1)^{\otimes k}$.  Hence we get $$(\bov,\  \bov_2)_{\scriptsize k} = (\bov, \ \sum_{i=1}^ne_{i+1,i}\bou_i)_{\scriptsize k}= ( \sum_{i=1}^n(e_{i,i+1}\bov,\ \bou_i)_{\scriptsize k}=0,$$ i.e, $(\bov,\ \bov')_{\scriptsize k}=(\bov,\ \bov_1)_{\scriptsize k}\ne 0$ and the Lemma is proved.
\end{pf}

\subsection{} \label{basis} From now on, we fix $k\in\bn$ and a partition $\xi=\xi_1\ge\cdots\ge\xi_{n+1}> 0$ of $k$. Let $k_i=\xi_i-\xi_{i+1}$ and set  $$\boe(\xi)=e_1^{\otimes k_1}\otimes (e_1\wedge e_2)^{\otimes k_2}\otimes\cdots\otimes (e_1\wedge e_2\wedge\cdots \wedge e_{n+1})^{\otimes k_{n+1}}\in  (V(\omega_1)^{\otimes k})^{\lie n^+}_\xi. $$ Let $S_k(\xi)$ be the isotropy subgroup of $S_k$  which fixing $\boe(\xi)$ and fix a set of distinct coset representatives $\sigma_\ell$, $1\le \ell \le d_\xi-1$ with $\sigma_1=\id$.
We note the following consequence of
Proposition \ref{schweyl} and the discussion in Section \ref{form1}.
 \begin{lem} The elements $\{\boe(\xi)\sigma_\ell: 1\le \ell\le d_\xi\}$ are a basis of $(V(\omega_1)^{\otimes k})^{\lie n^+}_\xi$. Moreover there exists a dual basis $\{\boe_\ell: 1\le \ell\le d_\xi\}$ of $(V(\omega_1)^{\otimes k})_\xi^{\lie n^+}$ satisfying \begin{equation}\label{dual}(\boe_j, \ \ \boe(\xi)\sigma_m)_{\scriptsize k}= 0,\ \ j\ne m,\qquad (\boe_m,\ \ \boe(\xi)\sigma_m)_{\scriptsize k}=1.\end{equation}

 \hfill\qedsymbol
 \end{lem}

\subsection{} We now turn our attention to representations of $\lie{sl}_{n+1}[t]$ and we begin with some general definitions. Given a  complex Lie algebra $\lie a$ and  an  indeterminate $t$, let  $\lie a[t]=\lie a\otimes \bc[t]$ be  the Lie algebra with commutator given by, $$[a\otimes f, b\otimes g]=[a,b]\otimes fg,\ \ a, b\in\lie a, \ \ f,g,\in\bc[t].$$  We identify without further comment the Lie algebra $\lie a$ with the subalgebra $\lie a\otimes 1$ of $\lie a[t]$. Clearly,  $\lie a[t]$ has a natural $\bz_+$--grading given by the powers of $t$ and this also induces a $\bz_+$--grading on $\bu(\lie a[t])$.
A $\bz_+$--graded bimodule for the pair $(\lie a[t], A_k)$ is a complex vector space $V =\oplus_{s\in\bz_+}V_s$ which admits a left action of $\lie a[t]$, a right action of $A_k$ and both actions are compatible with the grading, $$(\lie a\otimes t^r)V[s]\subset V[s+r],\qquad  V[s] A_k[r]\subset V[s+r],\ \ s,r\in\bz_+.$$


An elementary way to construct such modules is as follows: let $U_s$, $1\le s\le k$ be  $\lie a$--modules and    define on  $(U_1\otimes \cdots\otimes U_k\otimes A_k)$ the structure of  a  bimodule for the pair $(\lie a[t], A_k)$ by:
\begin{gather} (u_1\otimes \cdots\otimes u_k\otimes a)b=u_1\otimes \cdots\otimes u_k\otimes ab,\\
\label{actdef}(x\otimes t^r) (u_1\otimes \cdots\otimes u_k\otimes a)=\sum_{s=1}^k\left(\otimes_{j=1}^{s-1}u_j\right)\otimes xu_s\left(\otimes_{j=s+1}^{k}u_j\right)\otimes t_s^ra.\end{gather} The grading on $(U_1\otimes \cdots\otimes U_k\otimes A_k)$ is induced by the grading on $A_k$, i.e., for an integer $r$, the $r^{th}$--graded piece is $(U_1\otimes \cdots\otimes U_k)\otimes A_k[r]$.
These modules are clearly free as right $A_k$ (and so also as right $A_k^{S_k}$) modules. 

\subsection{} \label{bimod} Consider the special case when all the $U_i$'s are the same, say $U_i=U$ for $1\le i\le k$. The induced diagonal action of $S_k$ on $U^{\otimes k}\otimes A_k$
commutes with  the right action of $A_k$ and hence for each $\xi\in\cal P^k$, the isotypical component $(U^{\otimes k}\otimes  A_k)^{\xi}$ is a  $ A_k^{S_k}$--bimodule and we have an isomorphism of $A_k^{S_k}$--modules
 $$U^{\otimes k}\otimes A_k\cong \bigoplus_{\xi\in\cal P^k}(U^{\otimes k}\otimes  A_k)^{\xi},$$ If $\dim U<\infty$ then $(U^{\otimes k}\otimes  A_k)^{\xi}$ is a free $A_k^{S_k}$--module of finite rank. The following is easily checked by using the explicit definition of the action of $\lie a[t]$ given in \eqref{actdef}.
 \begin{lem} The space $(U^{\otimes k}\otimes A_k)^{S_k}$ is a $\bz_+$--graded    $(\lie a[t], A_k^{S_k})$--submodule of $U^{\otimes k}\otimes A_k$.\hfill\qedsymbol
 \end{lem}

\subsection{} Apply the preceding construction to the case when $\lie a=\lie{sl}_{n=1}$ and $U=V(\omega_1)$. Define   an $A_k$--valued bilinear form   $$\langle \ , \ \rangle:(V(\omega_1)^{\otimes k}\otimes A_k)\times (V(\omega_1)^{\otimes k}\otimes A_k)\to A_k,$$ by extending  the assignment,
$$\langle\bov \otimes f, \ \ \bov'\otimes g\rangle= (\bov, \bov')_{\scriptsize k} fg,$$ where $\bov,\bov'\in V(\omega_1)^{\otimes k}$,  and $f,g\in A_k$. Notice that the form is not $S_k$--invariant, however, we do have the following which is proved by a direct calculation.
\begin{lem}\label{inv2} For all $x\in\lie{sl}_{n+1}$, $f\in \bc[t]$, we have $$\langle (x\otimes f)(\bov\otimes g), \ \bov'\otimes g'\rangle\  = \ \langle(\bov\otimes g,\  (\Omega(x)\otimes f) (\bov'\otimes g')\rangle,$$ where $\bov,\bov'\in (V(\omega_1)^{\otimes k}$ and $g,g'\in A_k$. Further, $$\langle\ \ ,\  \bov\rangle: V(\omega_1)^{\otimes k}\otimes A_k\to A_k,$$ is a map of $A_k$--modules for all $\bov\in V(\omega_1)^{\otimes k}$.\hfill\qedsymbol
\end{lem}

\subsection{} For $k\in\bz_+$, set $$W(k)=\left(V(\omega_1)^{\otimes k}\otimes A_k\right)^{S_k}.$$
Since $W(k)$ is a locally finite--dimensional $\lie{sl}_{n+1}$--module we can write $$W(k)=\bigoplus_{\bor\in\bz^n} W(k)_\bor,\qquad W(k)_\xi=W(k)_\xi^{\lie n^+} \oplus(\lie n^-W(k)\cap W(k)_\xi). $$ Since the action of $\lie {sl}_{n+1}$ on $W(k)$ commutes with the action of
 $A_k^{S_k}$, we see that $W(k)_\bor$ and $W(k)_\xi^{\lie n^+}$ are graded free right $A_k^{S_k}$--module and in fact, we have \begin{equation}\label{globweyl}W(k)_\bor=  \left(V(\omega_1)^{\otimes k}_\bor\otimes A_k\right)^{S_k},\qquad \ W(k)_\xi^{\lie n^+}=  \left((V(\omega_1)^{\otimes k})_\xi^{\lie n^+}\otimes A_k\right)^{S_k}.\end{equation}

  \subsection{} We shall need the following result. As usual we denote by $\bu(\lie a)$ the universal enveloping algebra of the Lie algebra $\lie a$.

\begin{thm} \label{weylprop}Let $k\in\bn$.
\begin{enumerit}
\item[(i)] The element $v_{\omega_1}^{\otimes k}$ generates $W(k)$ as a module for  $\lie{sl}_{n+1}[t]$ and moreover, $$W(k)=\bu(\lie n^-[t])(v_{\omega_1}^{\otimes k}\otimes A_k^{S_k}).$$
    \item[(ii)] For $\xi\in\cal P^n$, the $A_k^{S_k}$--module $W(k)^{\lie n^+}_\xi$ is graded free and the number of elements in a  basis  with  grade $s$ is the coefficient of $q^s$ in the Kostka polynomial $\kappa_{1^k,\xi^{tr}}$.\end{enumerit}\hfill\qedsymbol
\end{thm}
Part (i) of the theorem was proved in  \cite{CPweyl}
     for $\lie{sl_2}$ and in \cite{FL} for $\lie{sl}_{n+1}$. Recall that   $I_k$ is the augmentation ideal in $A_k^{S_k}$ and consider  the graded $\lie{sl}_{n+1}[t]$--module, $(W(k)/I_kW(k)$. It was shown in \cite{CL} that the  subspace
     $(W(k)/I_kW(k)){\lie n^+}_\xi$ has a graded basis and the number of elements in the basis  of grade $s$ is the coefficient of $q^s$ in the Kostka polynomial $\kappa_{1^k,\xi^{tr}}$. Part (ii) follows since $W(k)^{\lie n+}_\xi$ is graded free as an $A_k^{S_k}$--module  and we have an isomorphism of graded spaces, $$W(k)^{\lie n+}_\xi/I_kW(k)^{\lie n^+}_\xi\cong (W(k)/I_kW(k))^{\lie n^+}_\xi.$$

\subsection{} The modules $W(k)$  are special examples  of a family of infinite--dimensional modules called
 global Weyl modules which are defined for arbitrary simple Lie algebras.
  For $\lie{sl}_{n+1}$ it is proved in \cite{BCGM} that the global Weyl modules
  can be realized as sitting inside a suitable space of invariants.
   However, except in the special case we consider in this paper, the global Weyl modules
    are strictly smaller than the space of invariants.

\section{The main result and Proof of Theorem 1}

\subsection{} As before we fix $k\in\bn$, $\xi=\xi_1\ge\xi_2\ge\cdots\xi_{n+1}>0$ and set $k_i=\xi_i-\xi_{i+1}$, $1\le i\le n+1$ and we freely use the notation of the earlier sections. Theorem 1 is clearly a consequence Theorem 2 and  the following result.

\begin{thm} The restriction of the map $\langle\ \ \ ,\ \ \boe(\xi)\rangle: V(\omega_1)^{\otimes k}\otimes k\to A_k$   gives an isomorphism of $A_k^{S_k}$--modules $W(k)^{\lie n^+}_\xi\to M_{k,\xi}$. Moreover,  for all $\sigma\in S_k$ we have \begin{gather}\label{saction} M_{k,\xi}\sigma=\langle W(k)^{\lie n^+}_\xi\ ,\ \boe(\xi)\sigma\rangle,\\ \label{isotyp}\langle W(k)^{\lie n^+}_\xi\ , \ (V(\omega_1)^{\otimes k})^{\lie n^+}_\xi\rangle\ \ \cong\ \ \bigoplus_{s=1}^{d_\xi}M_{k,\xi}\sigma_s\ \ \cong \ \ A_k^\xi,\end{gather} as $S_k$ and $A_k^{S_k}$--modules.
\end{thm}
\begin{pf} Since $\langle \ , \ \boe(\xi)\rangle: V(\omega_1)^k\otimes A_k\to A_k$ is a map of $A_k$--modules, it follows that the restriction to  $W(k)^{\lie n^+}_\xi$ is a map of $A_k^{S_k}$--modules. Using \eqref{globweyl} we have $$W(k)^{\lie n^+}_\xi = ((V(\omega_1)^k)^{\lie n^+}_\xi\otimes A_k)^{S_k},$$and hence we may write any non--zero element $\bov\in W(k)^{\lie n^+}_\xi$ as  $\bov=\sum_{s=1}^{d_\xi} \boe_s\otimes g_s,$ where  $g_s\in A_k$ and 
   $\{\boe_s: 1\le s\le d_\xi\}$ is the basis defined in
Section \ref{basis}.  By Lemma \ref{gp} we see that $g_1\ne 0$ and using \eqref{dual} we now see that $$\langle\bov,\ \boe(\xi)\rangle= g_1\ne 0,$$ which proves that the map $W(k)^{\lie n^+}_\xi\to A_k$ is injective.

  To prove that  $$\langle\ W(k)^{\lie n^+}_\xi,\ \boe(\xi)\rangle =M_{k,\xi},$$ we prove simultaneously that the elements $\bop(\bor)\in\langle\ W(k)^{\lie n^+}_\xi,\ \boe(\xi)\rangle $ and that they generate the $\langle\ W(k)^{\lie n^+}_\xi,\ \boe(\xi)\rangle$  as an $A_k^{S_k}$--module.
  By Theorem \ref{weylprop} we may write an element $\bov\in W(k)_\xi$ as  $$\bov=\sum_{s=1}^m y_s(v_{\omega_1}^{\otimes k}\otimes a_s),\ \  y_s\in\bu(\lie n^-[t]),\ \ a_s\in A_k^{S_k}.$$
   Since $x_{j,k}v_{\omega_1}=0,\ \ k\ne 1 \ \ ,$ it follows from the Poincare--Birkhoff--Witt theorem that we may assume that $y_s$ is a product of elements from the set $\{x_{j,1}\otimes t^r: 2\le j\le n+1,\ \ r\in\bz_+\}$. Moreover since $[x_{j,1}, x_{\ell,1}]=0$ for $2\le j,\ell\le n+1$ we can  write $$y_s=\prod_{j=2}^{n+1}(e_{j,1}\otimes t^{r_{j,1}})\cdots (e_{j,1}\otimes t^{r_{j,p_j}}), $$ for some $r_{j,p}\in\bz_+$, $2\le j\le n+1$  and $1\le p_j\le d_{\xi_j}$.
If $r_{j,p}=0$ for some $j$, i.e., if $y_s\in \lie n^-\bu(\lie n^-)[t]$, then since $$\langle\lie n^-\\bu(\lie n^-)[t]v_{\omega_1}^{\otimes k},\ \ boe(\xi)\rangle= \langle\bu(\lie n^-)[t]v_{\omega_1}^{\otimes k},\ \ \lie n^+boe(\xi)\rangle =0,$$ we see that
   \begin{equation}\label{zero} (y_sv_{\omega_1}^{\otimes k}\ ,  \ \boe(\xi))= 0. \end{equation} If  $r_{j,p}>0$ for all $j$ then,
     $$\langle y_s v_{\omega_1}^{\otimes k} \ , \  \boe(\xi)\rangle =\langle v_{\omega_1}^{\otimes k}\ , \ \prod_{j=2}^{n+1}(e_{1,j}\otimes t^{r_{j,1}})\cdots (e_{1,j}\otimes t^{r_{j,p_j}})\boe(\xi)\ \rangle = \bop(\bor),$$ where 
    a tedious but straightforward calculation shows that \begin{equation}\label{equal} \prod_{j=2}^{n+1}(e_{1,j}\otimes t^{r_{j,1}})\cdots (e_{1,j}\otimes t^{r_{j,p_j}})\boe(\xi)\  = v_{\omega_1}^{\otimes k}\otimes \bop(\bor).\end{equation}
For instance in the case of $\lie{sl}_2$, one sees that $$ (e_{1,2}\otimes t^{r_1})\otimes\cdots\otimes (e_{1,2}\otimes t^{r_\ell})\left((e_1^{k-2\ell}\otimes (e_1\wedge e_2)^{\otimes \ell}\right)=\\ e_1^{k-2\ell}\otimes(e_{1,2}\otimes t^{r_1})\otimes\cdots\otimes (e_{1,2}\otimes t^{r_\ell})\left (e_1\wedge e_2)\right)^{\otimes \ell}.$$
    Now use the fact that $$e_{1,2}\otimes t^r(e_1\otimes e_2-e_2\otimes e_1)=(e_1\otimes e_1\otimes t_2^{r}-e_1\otimes e_1\otimes t_1^r=(e_1\otimes e_1)\otimes (t_2^r-t_1^r),$$ to complete the calculation.

A consequence of equations \eqref{zero} and \eqref{equal} is the following: if $y_sv_{\omega_1}^{\otimes k}\ne 0$, then    $y_sv_{\omega_1}^{\otimes k}$ has non--zero projection onto $W(k)^{\lie n^+}_\xi$. We have already noted that $\langle \lie n^-W(k),\boe(\xi\rangle=0$ and hence we have now proved both that $\bop(\bor)$ is in the image of the map $\langle \ \ ,\  \boe(\xi)\rangle$ and that the image is generated as an $A_k^{S_k}$--module by these elements.

\medskip

To prove \eqref{saction},  let  $\boa\in M_{k,\xi}$ and   choose $\bov=\sum_{s=1}^\ell v_s\otimes f_s\in W(k)^{\lie n^+}_\xi$ such that $$\boa=\langle \bov\ ,\ \boe(\xi)\rangle=\sum_{s=1}^\ell( v_s,\ \boe(\xi))_{\scriptsize k}f_s.$$ For $\sigma\in S_k$, we get by using the $S_k$--invariance of $(\ ,\ )_{\scriptsize k}$ that, \begin{gather*}\boa\sigma=\sum_{s=1}^\ell( v_s\sigma,\ \boe(\xi)\sigma)_{\scriptsize k}f_s\sigma= \langle\sum_{s=1}^\ell v_s\sigma\otimes f_s\sigma,
\boe(\xi)\sigma\rangle\\ =\langle(\sum_{s=1}^\ell v_s\otimes f_s)\sigma,
\boe(\xi)\sigma\rangle=\langle\bov\sigma, \boe(\xi)\sigma\rangle=\langle\bov\ ,\ \boe(\xi)\sigma\rangle,\end{gather*}where the last equality is  consequence of the fact that elements of $W(k)$ are $S_k$--invariant. This proves that
  $$M_{k,\xi}\sigma=\langle W(k)^{\lie n^+}_\xi,\boe(\xi)\sigma\rangle,\ \ \sigma\in S_k,$$ and hence also that we have an isomorphism of $S_k$ and $A_k^{S_k}$ modules,$$\langle W(k)^{\lie n^+}_\xi\ , \ (V(\omega_1)^{\otimes k})^{\lie n^+}_\xi \rangle \cong \sum_{\sigma\in S_k}M_{k,\xi}\sigma.$$

We now prove that $$\sum_{\sigma\in S_k}M_{k,\xi}\sigma=A^\xi.$$ To prove that the left hand side is contained in the right hand side, it suffices to prove that $M_{k,\xi}\subset A^\xi$. But this is clear since as we have observed before any  $\boa\in M_{k,\xi}$ may be written as,
 $$\boa=\langle \bov\ ,\ \boe(\xi)\rangle=\sum_{s=1}^\ell( v_s,\ \boe(\xi))_{\scriptsize k}f_s,$$ where $$\bov\in W(k)^{\lie n^+}_\xi=\left (V(\omega_1)^{\otimes A_k})^{\lie n^+}_\xi\otimes A_k\right)^{S_k}= \left (V(\omega_1)^{\otimes k})^{\lie n^+}_\xi\otimes A^\xi_k\right)^{S_k}.$$ This means we may assume that $f_s\in A^\xi$ for all $1\le s\le \ell$ which proves that $\boa\in A^\xi$.

For the reverse inclusion, suppose that $a\in A_k^\xi$ and assume that it generates an irreducible $S_k$--submodule $N$  and let $a_1,\cdots, a_{d_\xi}$ be a basis of $N$ where $a_1=a$. Then there exists a non--zero element  $$\bov\in (V(\omega_1)^{\otimes k})^{\lie n^+}_\xi\otimes N)^{S_k}\subset W(k)^{\lie n^+}_\xi,$$ and by Lemma \ref{gp} we can  write $ \bov=\sum_{s=1}^{d_{\xi}}\bov_s\otimes a_s$ with  $ \bov_s\ne 0$. In fact, the elements $\{\bov_s:1\le s\le d_\xi\}$ are a basis for $(V(\omega_1)^{\otimes k})^{\lie n^+}_\xi$ and by Lemma \ref{nondeg} we may choose a dual basis  $\{\bov'_s:1\le s\le d_\xi\}$ and we find now that $$\langle\bov,\ \bov_1'\rangle =a_1=a,$$ which proves that $A^\xi$ is contained in $M_\xi$.

To complete the proof of the theorem, we must prove that $M_{k,\xi}\sigma$ is a summand of $A^\xi$ as an $A_k^{S_k}$--module. We have so far proved that $M_{k,\xi}\sigma$ is free of rank $d_\xi$ and also that $$A^\xi=\sum_{s=1}^{d_\xi} M_{k,\xi}\sigma_s.$$ Since $A^\xi$ is a free module of rank $d_\xi^2$, it follows that the the sum is direct and the proof of the theorem is complete.

\end{pf}

\end{document}